\documentclass[11pt]{article}   
\usepackage{amssymb,amscd,latexsym}   
\usepackage{amsmath}
\usepackage{amsthm}
\input xypic
\newcommand{\llar}{-\kern-5pt-\kern-5pt\longrightarrow}

\newtheorem{Theorem}{Theorem}[section]
\newtheorem{Lemma}[Theorem]{Lemma}
\newtheorem{Corollary}[Theorem]{Corollary}
\newtheorem{Proposition}[Theorem]{Proposition}
\newtheorem{Remark}[Theorem]{Remark}
\newtheorem{Example}[Theorem]{Example}

\newtheorem{Definition}[Theorem]{Definition}

\def\sqr#1#2{{\vcenter{\hrule height.#2pt
        \hbox{\vrule width.#2pt height#1pt \kern#1pt
            \vrule width.#2pt}
        \hrule height.#2pt}}}
\def\phi{\varphi}
\def\demo{\noindent{\bf Proof. }}
\def\square{\mathchoice\sqr64\sqr64\sqr{4}3\sqr{3}3}
\def\qed{\hspace*{\fill} $\square$}

\def\xx{{\bf x}}
\def\yy{{\bf y}}

\def\XX{{\bf X}}

\def\Q{{\mathbb Q}}
\def\R{{\mathbb R}}

\def\Z{{\mathbb Z}}
\def\NN{\mathbb N}

\def\pp{{\mathbb P}}

\begin{document}
\begin{center}
{\Large{\bf\sc Cremona maps defined by monomials}}\footnotetext{Mathematics Subject Classification 2000 (MSC2000).
Primary 14E05, 14E07, 13A30, 52B20;
Secondary 05B20, 05C50, 11C20, 13B22, 15A36,  14E99.}

\vspace{0.3in}

{\large\sc Barbara Costa}\footnote[1]{Under a CNPq Doctoral scholarship},
 and {\large\sc Aron  Simis}\footnote[2]{Partially
supported by a CNPq grant.}

\end{center}


\bigskip

\begin{abstract}

Cremona maps defined by monomials of degree $2$ are thoroughly analyzed
and classified via integer
arithmetic and graph combinatorics. In particular, the structure of the inverse
map to such a monomial Cremona map is made very explicit as is the degree
of its monomial defining coordinates.
As a special case, one proves that any monomial Cremona map of degree $2$
has inverse of degree $2$ if and only if it is an involution up to permutation
in the source and in the target. This statement is subsumed
in a recent result of L. Pirio and F. Russo, but the proof is entirely different and
holds in all characteristics.
One unveils a close relationship binding together the normality
of a monomial ideal, monomial Cremona maps and Hilbert bases of polyhedral cones.
The latter suggests that facets of monomial Cremona theory may be {\bf NP}-hard.

\end{abstract}

\section*{Introduction}

The expression ``birational combinatorics'' has first been used in \cite{birational-linear}
to mean the theory of
characteristic-free rational maps $\pp^{n-1}\dasharrow \pp^{m-1}$ defined by
monomials, along with natural criteria for such maps to be birational onto their
image varieties. As emphasized in [loc.cit.], the theory and the criteria were designed to
reflect the specificity of the relevant combinatorial data, thus retracting from
the classical theory of Cremona transformations in characteristic zero.

The challenge remained, as said in \cite{birational-linear}, as to how one would proceed to find the inverse map to
a monomial birational map by
a purely combinatorial method.
Since the inverse has to be given by monomials as well, as shown in {\em ibid.}, the question made sense.
This question was eventually solved in \cite{CremonaMexico} a few years later.

One of the peculiarities of the theory is that even if the given monomials are squarefree to start with,
the inverse map is pretty generally defined by non-squarefree monomials.
This makes classification in high degrees, if not the structure of the monomial Cremona group itself, a tall order.

In this paper, we continue along these steps by tackling the following questions:

\begin{itemize}
\item Classification of monomial Cremona maps of degree $2$ in any number of variables
\item The structure of the inverse map to a monomial Cremona map of degree $2$
\item The role of Hilbert bases in monomial Cremona maps of arbitrary degrees.
\end{itemize}

The first two questions hinge on a certain normal form for the so called {\em log-matrix} of a set
of monomials of degree $2$. Quite generally, this matrix plays an essential role  in the
combinatorial criteria obtained in the previous references.
Thus, we proceed {\em ab initio} by giving a reasonable unique normal form of a monomial Cremona map
of degree $2$ based on the structure of the graph whose edges correspond to the defining monomials.
The overall motivation has been to obtain the explicit format of the inverse map and the associated numerical
invariants in terms of the nature of the corresponding graph.
This is accomplished in our main result (Theorem~\ref{degree_of_inverse}).

We have also been driven by the a question  related to
the monomial Cremona maps of degree $2$ whose inverse has also degree $2$.
We have been able to characterize these Cremona maps in terms
of involutions and the nature of the corresponding graph.
We were originally motivated by a conversation with F. Russo about this matter.
His joint results with L. Pirio deal with general such Cremona maps and give an important connection
with the theory of Jordan algebras (see \cite[5.3, 5.13]{PiRu}).
They prove, among other things, that any Cremona map of degree $2$ whose inverse is also
of degree $2$ is an involution up to a projective change of coordinates.
Though a lot simpler our present result in the monomial case is characteristic free.
It would be nice to prove the above result of Pirio--Russo in all characteristics.

Though the criteria themselves have an expected simplicity and afford effective computation (see \cite[Section 4]{CremonaMexico})
-- although facets of the computation are close to {\bf NP}-hard problems -- the practical use in theoretical classification is by no means obvious,
often requiring quite a bit of ingenuity.
Having traded geometrical tools by integer linear algebra one pays a price in that further precision
has to be exercised. This sort of toil will be found throughout the arguments of the main results.

\medskip

We now describe more closely the contents of each section

The first section is about the background terminology and a review of previous results of
the birational combinatorics repository. This will hopefully help increase familiarity with
the language in which the subsequent results are stated and also make the needed acquaintance with the
fundamental facts to be used throughout.
The one main criterion drawn upon is carefully stated as Theorem~\ref{teorema_fundamental}.
This result explains, in particular, the nature of the inverse map of a monomial Cremona map in terms of the
corresponding log-matrix, here called perhaps appropriately {\em Cremona inverse matrix}.
The source log-matrix and its Cremona inverse are both stochastic and the theorem gives a nice
equation connecting the corresponding stochastic numbers -- these numbers coincide with
the respective degrees in the traditional Cremona terminology. As in the classical case
of arbitrary Cremona maps in characteristic zero, the degrees are related through properties
of the base locus. However, in the monomial case, this relation is really an equation whose
terms are effectively computable  in terms of optimization methods, however intractable
as they may be from the viewpoint of computational complexity.

\medskip

The second section is about Cremona maps defined by monomials of degree $2$ and constitutes
the core of the paper.
Though simple to grasp in terms of the corresponding graph, such a map keeps a couple of
hidden marvels, such as the precise format and degree of its Cremona inverse.
Moreover, it touches some pertinent theoretic aspects along the way, such as to when
the base ideal of a Cremona map is of linear type.
Cremona maps whose base ideals are of linear type have been considered in recent work
(see \cite{RuSi}, \cite{bir2003} and \cite{AHA}).
Our first result in this section (Proposition~\ref{linear_type_of_inverse}) gives a complete characterization of when
the base ideal of the Cremona inverse to a monomial Cremona generated in degree $2$
is of linear type. According to a later terminology introduced in our development, there are
just ``a few'' of these.

The main result of the section is Theorem~\ref{degree_of_inverse} which determines the degree
of the inverse to a monomial Cremona generated in degree $2$ in terms of the corresponding
graph. The proof is elaborate and long, however it has the advantage of producing along its
way the precise format of the Cremona inverse matrix.
Moreover, as relevant by-products we obtain the precise format of the so-called {\em inversion
factor} and an easy criterion in terms of the associated graph as to when the inverse map
is defined by squarefree monomials.

There are many consequences of this theorem to classification, of which the main ones are
Proposition~\ref{involution_deg_2}, Corollary~\ref{graph_types}
and Corollary~\ref{graph_of_apocryphal}.
The terminology becomes slightly technical to be shortly explained in this introduction, so we
refer to the appropriate parts in the section.
Here, classification means to uncover certain classes of specific behavior among all Cremona maps
defined by monomials of degree $2$. It is not clear what is the impact of this sort of
classification on the structure of the monomial Cremona group (but see \cite{Barbara} for a couple of hints).

\medskip

The third section is entirely devoted to developing a problem suggested by R. Villarreal.
The outcome is a curious relation between the notion of a normal ideal (in the monomial case)
and Hilbert bases in the sense of combinatorics.
The required technical parts are  kept to  a minimum and the bridging is discussed as much as possible,
hopefully without jamming the overall reading.
The main results of this section are Theorem~\ref{main_normal} and Theorem~\ref{normal2Hilbert2Cremona}.
The second of these theorems shows that, given a monomial birational map $\pp^n\dasharrow \pp^m$
(onto the image) whose base ideal is normal,
there exists a coordinate projection $\pp^m\dasharrow \pp^n$ such that the composite is a Cremona map.
The proof draws upon Hilbert base technique and s couple of arithmetical lemmas established in previous
work on birational combinatorics (\cite{SiVi}).

\section{Terminology and basic combinatorial criterion}

In this short section we state the setup and the terminology of monomial rational maps,
as well as a basic integer arithmetic
formulation of monomial Cremona transformations.

Let $k$ be a field and let $R=k[\xx]=k[x_1,\ldots,x_n]$ a polynomial ring over $k$.
Though not strictly needed for the preliminaries, the theory will be meaningless for $n=1$. Thus, we assume once for all
and without further say, that $n\geq 2$.

Given $\alpha=(a_1,\ldots,a_n)\in {\mathbb N}^n$, write ${\xx}^{\alpha}:= x_1^{a_1}\cdots x_{n}^{a_{n}}$
for the associated monomial.
We will be concerned with a finite set $V=\{v_1,\ldots,v_q\}\subset \mathbb{N}^n$ of distinct vectors and
the corresponding {\em log\/}-set of monomials $F^V=\{{\xx}^{v_1},\ldots, {\xx}^{v_q}\}\subset R$.
The following basic restrictions will be assumed throughout:

\begin{itemize}
\item For every $i\in \{1,\ldots,n\}$ there is a $j\in\{1,\ldots,q\}$ such that $v_{ij}=0$
\item For every $i\in \{1,\ldots,n\}$ there is at least one $j\in\{1,\ldots,q\}$ such that $v_{ij}\neq 0$
\end{itemize}

For convenience of cross reference, we will call this set of assumptions the {\em canonical restrictions}.

The first requirement is not so intuitive, but it becomes clear in terms of the corresponding
$F^V$. It means that the monomials in $F^V$ have no non-trivial common factor.

The second requirement can always be achieved by simply contracting to a subset of a suitable
coordinate ${\mathbb N}^{n-1}\subset {\mathbb N}^n$.
(This requirement is analogous to the idea of a subvariety that is not a cone in the sense of algebraic geometry.)

\begin{Definition}\rm The {\em log-matrix} $A_V$ of $V$ as above is the integer matrix whose columns
are the vectors in $V$.
\end{Definition}

By extension, one calls $A_V$ the log-matrix of the corresponding set $F=F^V$ of monomials.
Accordingly, we will throughout use the notation $A_F$.

There is of course a slight instability in this terminology, as the matrix depends
on the order of the variables and on the order of the monomials in $F$ -- the same sort of imprecision
that one faces when talking about {\em the} Jacobian matrix of a set of polynomials.
This instability is most of the times harmless, although some care has to be exercised, specially
in the statements of results and in arguments that involve several sets of monomials.

\medskip

We will be exclusively dealing with stochastic sets of vectors. More precisely,
$V$ is a $d$-{\em stochastic} set in the sense that $|v_1|=\cdots =|v_q|=d$,
for some fixed integer $d\geq 1$.

To see how the log-matrix of a $d$-stochastic set of vectors comes about in the subject, recall
that an extension $D'\subset D$ of integral domains is said to be birational if it is an
equality at the level of the respective fields of fractions. Let
$F=F^V\subset R=k[\xx]$ be the set of monomials associated to a $d$-stochastic set $V\subset \mathbb{N}^n$.
Write ${\bf x}_d$ for the set of {\it all\/} monomials of degree $d$ in $R$. Then
$k[{\bf x}_d]$ is the $d${\it th} Veronese subring $R^{(d)}$ of $R$. The overall aim
is to understand the birationality of the ring extension $K[F]\subset R^{(d)}$ -- this can be translated
into geometry as the birationality of $\pp^{n-1}$ onto the image of the rational map
$\pp^{n-1}\dasharrow \pp^{q-1}$  whose defining coordinates are $({\xx}^{v_1}:\cdots: {\xx}^{v_q})$.

The fundamental transition from algebra/geometry to integer arithmetic is processed through
the following simple result.

\begin{Lemma}\label{birational_principle}{\rm (\cite[Proposition 1.2]{SiVi})}
Let $V\subset \mathbb{N}^n$ be a $d$-stochastic set satisfying the canonical restrictions,
with $d\geq 1$.
Then $k[{F^V}]\subset R^{(d)}$ is a birational extension if and only if
the ideal of $\Z$ generated by the $n\times n$-minors of $A_F$ is generated by $d$.
\end{Lemma}

Our goal here is more restricted in that we assume the case $q=n$.
If the map is birational, one calls it a {\em Cremona map} or {\em Cremona transformation}, a venerable classical object.
We are then dealing with monomial maps  whose corresponding log-matrix is an $n\times n$ $d$-stochastic matrix of
determinant $\pm d$.

The sign above is not relevant as one can always permute two columns to achieve a positive determinant, a harmless
operation corresponding to a transposition in the set of the monomials in $F$.
Actually, often the results may depend on allowing a permutation of the set of variables and the set
of monomials in $F$.
This means that, in the monomial Cremona group in $n$ variables, we will not distinguish between an element
$F$ of this group and the composite $PFQ$, where $P,Q$ are arbitrary {\em Cremona permutations} -- i.e., those
whose corresponding log-matrices are permutation matrices.
In particular, it is clear that across such compositions we are not changing the degree of the defining monomials.

\medskip

We will also make heavy use throughout of the main result in \cite{CremonaMexico}
giving the integer arithmetic counterpart of a monomial Cremona situation:

\begin{Theorem}{\rm (\cite[Theorem 2.2]{CremonaMexico})}\label{teorema_fundamental}
Let $V=\{v_1,\ldots,v_n\}\subset \mathbb{N}^n$ stand for a $d$-stochastic set, with $d\geq 1$, satisfying the
canonical restrictions.
Suppose that the determinant of the associated log-matrix $A_V$ is $\pm d$.
 Then there exists a unique set $W=\{w_1,\ldots,w_n\}\subset \mathbb{N}^n$ and a unique vector
 $\gamma\in \mathbb{N}^n$  such that{\rm :}
\begin{enumerate}
\item[\rm (a)] $A_V\cdot A_W=\Gamma+I_n$, where $\Gamma=[\underbrace{\gamma|\cdots |\gamma}_n]\,${\rm ;}
\item[\rm (b)] $W$ is $\delta$-stochastic and satisfies the canonical restrictions, and $\det(A_W)=\pm\delta$, where $\delta=\frac{|\gamma|+1}{d}$.
\end{enumerate}
\end{Theorem}
If $F$ and $G$ are the respective log-sets of monomials, in order  to emphasize the combinatorial stage background,
we will also write the fundamental
equation of the theorem in the form $A_F\cdot A_G=\Gamma+I_n$, and refer to it as the {\em inversion equation}
of $F$, while $A_G$ is  referred to as the {\em Cremona inverse matrix} of $A_F$.

The theorem shows in particular that the inverse of a monomial Cremona map is also monomial,
thus reproving the fact that the subset of the entire Cremona group whose elements are the
monomial Cremona maps is a subgroup.
Yet more important is that its proof gives explicitly the Cremona inverse matrix.
Moreover, the result ties both Cremona matrices by means of an integer vector yielding the
proportionality monomial factor responsible for the composition of the maps being
the identity.
We will have more to say about this factor in later sections.
Whenever needed we refer to this vector (respectively, factor) as the {\em inversion vector}
(respectively, the {\em inversion factor}) of $F$.

\section{The inverse of a Cremona map of degree $2$}

Parts of this work will concern the special case $d=2$.
In this situation, there is a natural graph $G_F$ associated to $F$, whose set of vertices
is in bijection with $\{x_1,\ldots,x_n\}$ and whose set of edges corresponds to the set $F$
in the obvious way.
Note that $G_F$ may have loops, corresponding to pure powers of order $2$ among the monomials
in $F$.
Also, the log-matrix $A_F$ is exactly the incidence matrix of the graph $G_F$.

\medskip

The following notion has been introduced in \cite{birational-linear}.
A set $V\subset \mathbb{N}^n$ is {\em non-cohesive} if, up to permutation of either rows or columns,
the log-matrix is block-diagonal, i.e., of the form
$$A_V=\left(
\begin{array}{c@{\quad\vrule\quad}c}
C &  \\
\multispan2\hrulefill\\
  & D
\end{array}
\right),$$
where $C,D$ are log-matrices of suitable sizes and the empty slots have only zero entries.
We say that $V$ is {\em cohesive} in the opposite case.
By extension, if $V$ is cohesive we often say that the corresponding set $F^V$ of monomials is cohesive.
Particularly, if $V$ is $2$-stochastic, then cohesiveness of $V$ translates into connectedness of the associated graph $G_F$.

We observe that, quite generally, a stochastic set $V$ such that $F^V$ defines a birational map onto its image
is necessarily cohesive (\cite[Lemma 4.1]{birational-linear}).

With this notion available, the basic dictionary above can be stretched to accommodate a graph-theoretic characterization
and an ideal-theoretic property.
Recall that an ideal $I\subset R$ is said to be of {\em linear type} if the natural surjective algebra homomorphism
from the symmetric algebra of $I$ to its Rees algebra is injective (i.e., an isomorphism).

\begin{Theorem}{\rm (\cite[Proposition 5.1]{birational-linear})}\label{cremona_grafo}
Let $F\subset k[x_1,\ldots, x_n]$ be a cohesive finite set of monomials of degree
$2$ having no non-trivial common factor and let $G_F$ denote the corresponding graph as above.
 The following conditions are equivalent:
\begin{enumerate}
\item[{\rm (i)}] $\det A_F\neq 0$
\item[{\rm (ii)}] $F$ defines a Cremona transformation of $\pp^{n-1}$
\item[{\rm (iii)}] Either
\begin{itemize}
\item[{\rm (a)}] $G_F$ has no loops and contains a unique circuit
 that is necessarily of odd length,
\end{itemize}
\vskip-8pt
 or else
\begin{itemize}
\item[{\rm (b)}] $G_F$ is a tree plus exactly one loop.
\end{itemize}
\item[{\rm (iv)}] The ideal $(F)\subset k[x_1,\ldots,x_n]$ is of linear type.
\end{enumerate}
\end{Theorem}

It is convenient to call the set $F$ a {\em Cremona set in degree} $2$ in the case it defines a Cremona map.
Likewise, the corresponding log-matrix can be called a {\em Cremona matrix in degree} $2$.

\subsection{The basic log-matrix in degree $2$}

In this part we introduce a normal form of a Cremona matrix in degree $2$.
Although  the inverse will only infrequently be in degree $2$,
the present normal form will shed light on the form of the corresponding
log-matrix.
Note that in degree $2$ the log-matrix is just the usual incidence matrix of the corresponding graph.

\begin{Lemma}\label{format_logmatrix}
Let $F$ be a Cremona set in degree $2$.
Then, up to permutations of the variables and the monomials of $F$, its log-matrix can be written
in the form
$${\cal N}_F = \left(
\begin{array}{c@{\quad\vrule\quad}c@{\quad\vrule\quad}c@{\quad\vrule\quad}c@{\quad\vrule\quad}c} 
{{\cal N}_F}_{_r} & M_1 & &  &  \\
\multispan5\hrulefill\\
  & I_{s_1} & M_2 &  &  \\
\multispan5\hrulefill\\
  &  & I_{s_2} &  &   \\
\multispan5\hrulefill\\
 & & & \ddots & M_p  \\
\multispan5\hrulefill\\
& & & & I_{s_p}  
\end{array}
\right),$$
where:
\begin{enumerate}
\item[{\rm (i)}]  ${{\cal N}_F}_{_r}$ denotes either the incidence matrix of the unique odd circuit
of length $r$, written in the form
$$\left(
\begin{array}{ccccc}
1 & 0 & & 0 & 1\\
1 & 1 & & 0 & 0\\
0 & 1 & & 0 & 0\\
 &  & \ddots & &\\
0 & 0 & &   1 & 0\\
0 & 0 & & 1 & 1
\end{array}
\right)
$$
 or the $1\times 1$ matrix $(2)$ corresponding to the unique loop
\item[{\rm (ii)}]   $I_{s_j}$ denotes the identity matrix of size $s_j$,
where $j$ runs through the ordered list of the neighborhoods of the unique circuit
and $s_j$ denotes the number of edges in the $j$th neighborhood of the circuit, with $j=1,\ldots,p$
\item[{\rm (iii)}]  $M_j$ is a certain matrix with $0,1$
entries distributed in such a way as to have exactly one $1$ on every column
\item[{\rm (iv)}]  The empty slots are filled with zeros.
\end{enumerate}
\end{Lemma}
\demo The proof is essentially contained in the itemized details of the statement.
For this one considers the graph $G_F$ as a rooted tree, where the root is either the unique circuit or else the unique loop, according to the case,
while the successive neighborhoods of the circuit (or loop) form the branching set of the tree. The disposition on the matrix
(left to right) obeys the order of the neighborhoods from the root.
Finally, along every single neighborhood we order the vertices  in an arbitrary way.
\qed

\begin{Definition}\rm The above matrix will be called a {\em normal form} of the Cremona matrix in degree $2$.
\end{Definition}

Clearly, this normal form is not unique as it depends on arbitrary ordering of vertices on each neighborhood of the root.
Actually, this degree of arbitrariness may turn out to be useful in an argument.

\medskip

Recall that, given a rational map $\pp^{n-1}\dasharrow \pp^{n-1}$, its {\em degree} is the common degree of its coordinate
forms, provided these forms have no proper common factor.
This is degree is not to be confused with the degree of the map in the geometric or field-theoretic sense.
Since we are dealing mostly with Cremona maps, the degree in the latter sense is always $1$, so there is no room for
misunderstanding.

In \cite[Proposition 5.5]{birational-linear} a classification has been given of Cremona maps of degree $3$ in $5$ variables.
By complementarity, this  was based on the corresponding classification in degree $2$.
In this section we intend to extend this classification to degree $2$ in any number of variables.

For this we will establish some preliminaries concerning the nature of the inverse map.

\subsection{When is the inverse map of linear type?}

We refer to the notation introduced in the previous subsection and, particularly, the notation employed
in the statement of Lemma~\ref{format_logmatrix}.

Recall that the {\em edge graph} of a graph $G$ is the simple graph whose vertex set is
the set of edges of $G$ and two such vertices
are adjacent if and only if the original edges meet (see \cite[Definition~6.6.1]{VillaBook}).
Observe that $G$ is connected if and only if its edge graph is connected.
The other notion we need is that of the {\em diameter} of a simple graph, defined to
be the largest distance between any two vertices of the graph, where the distance
between two vertices is the minimum number of edges connecting them.

We will need the following simple result.

\begin{Lemma}\label{diameter_of_cremona}
Let $F$ denote a Cremona set and let $G=G_F$ stand for the corresponding graph.
Then the diameter of the edge-graph of $G$ is bounded below by $\frac{r-1}{2}+p$,
where $r$ is the length of the unique circuit of $G$ and $p$ is the number of
ordered neighborhoods of the circuit.
\end{Lemma}
\demo Actually, the argument will show that the claimed bound is always attained, but
we will have no need for the equality.

First it is quite elementary that a circuit is edge-dual, i.e., it edge-graph is again a circuit of
the same length.
Clearly, then the diameter of the edge-graph of a circuit of length $r$ is $(r-1)/2$.

Write $\mathcal{E}(G), \mathcal{E}(C)$ for the edge-graph of $G$ and $C$, respectively, where $C$ denotes
the unique circuit of $G$ (note that, in case $C$ degenerates to a loop, $\mathcal{E}(C)=\emptyset$).

Let $v(p)$ be a vertex on the $p$th -- i.e., the last nonempty -- neighborhood  of the circuit in $G$
and let $v(p-1), v(p-2),\ldots,v(0)$ denote the unique set of  vertices belonging, respectively,
to the neighborhoods of $C$ of order $p-1, p-2,\ldots, 0$,  and such that $v(i)$ and $v(i+1)$
are adjacent for $i=0,\ldots,p-1$.

Then, in $\mathcal{E}(G)$ the distance between the ``edges'' $v(p)v(p-1)$ and $v(1)v(0)$ is exactly $p$
as there is no shorter path due to the ordered structure of neighborhoods of $C$ on $G$.
To still connect to an arbitrary vertex (``edge'') of $\mathcal{E}(C)$ will require at least $(r-1)/2$ additional
vertices (``edges'').
Therefore, the diameter of $\mathcal{E}(G)$  is at least $\frac{r-1}{2}+p$.
\qed

\medskip

The foregoing affords a graph-theoretic characterization of a Cremona set in degree $2$ whose
inverse is defined by an ideal of linear type.

\begin{Proposition}\label{linear_type_of_inverse}
Let $F\subset k[\xx]$ denote a Cremona set in degree $2$
satisfying the canonical restrictions and let $G=G_F$ stand for the corresponding graph, with unique circuit of
length $r$, possibly degenerating to a loop.
Let $s$ stand for the number of vertices
of $G$ off the circuit whose vertex degree is $\geq 2$.
Let $F^{-1}\subset k[\xx]$ denote the unique Cremona inverse set
satisfying the canonical restrictions,
The following conditions are equivalent:
\begin{enumerate}
\item[{\rm (a)}] The ideal $(F^{-1})\subset k[\xx]$ is of linear type
\item[{\rm (b)}] The ideal $(F)\subset k[\xx]$ is linearly presented
\item[{\rm (c)}] The edge-graph of $G$ has diameter $\leq 2$
\item[{\rm (d)}] The root-neighborhood structure of $G$ is one of the following:
\[\left\{
\begin{array}{ll}
r=1, & \mbox{$s_j=0$ for $j\geq 3$, and $s\leq 1$}\\
r=3 & \mbox{and $s_j=0$ for $j\geq 2$}\\
r=5 & \mbox{and $s_j=0$ for $j\geq 1$}.
\end{array}
\right.
\]
\end{enumerate}
\end{Proposition}
\demo  (a) $\Leftrightarrow$  (b)
This is because the respective Rees algebras of $(F)$ and $(F^{-1})$ have the the same ideal of relations
(\cite{bir2003}).
More precisely, taking a clone $\yy$  of the $\xx$-variables, the defining ideal of the Rees algebra of $(F)\subset k[\xx]$
on $k[\xx,\yy]$ is the same as the defining ideal of the Rees algebra of $(F^{-1})\subset k[\yy]$ on $k[\yy,\xx]$.
In particular, if there is some homogeneous syzygy of $(F)$ that is not generated by the linear ones it produces a relation in degree $2$ or higher
that is not generated by the degree $1$ relations of $(F^{-1})$, hence the latter could not be of linear type.
The converse is similar.

(b) $\Leftrightarrow$ (c)
This is pretty general (see, e.g., \cite[Lemma 5.16]{BGS}).

(d) $\Rightarrow$ (c)
This is by a direct inspection.
Namely, as previously observed, the pentagon is edge-dual, i.e., its edge-graph is also a pentagon, hence
has diameter $2$.

Next, a triangle with empty second neighborhood has an edge-graph which is
``sufficiently'' triangulated. Indeed, the triangle is also edge-dual, while any edge
issuing from the triangle adds another triangle in the edge-graph sharing a common edge
with the base triangle.
Thus, it is apparent that the diameter of the edge-graph is $2$.

Finally, for the case of the loop note again that the {\em star} subgraph rooted in the loop of $G$
gives rise to a complete subgraph in $\mathcal{E}(G)$ and so does the star subgraph
rooted on the unique  non-looped vertex $v$ of $G$ of degree $\geq 2$.
Moreover, these two complete subgraphs of $\mathcal{E}(G)$ meet on the unique vertex (``edge'')
joining $v$ to the loop vertex.
Therefore, the the diameter of $\mathcal{E}(G)$ is at most $2$.

(c) $\Rightarrow$ (d)
By Lemma~\ref{diameter_of_cremona}, we must have $r-1+2p\leq 4$.
This immediately forces the three stated  alternatives -- note that in the case of the loop,
$G$ has at most one non-loop vertex of degree $\geq 2$. Indeed, otherwise, given distinct vertices
$v(1), u(1)$ on the first neighborhood of the loop and two (necessarily distinct) vertices $v(2), u(2)$
adjacent to the first two, respectively, then the distance in $\mathcal{E}(G)$ between the ``edges''
$v(1)v(2)$ and $u(1)u(2)$ is $3$.
\qed

\subsection{The degree of the inverse map}

The next result draws on the normal form in order to tell  about the form of the inverse Cremona matrix - note it will
quite generally be a matrix in higher degree.
If $F$ is the set of monomials of a Cremona map, we write $F^{-1}$ for the set of monomials defining the inverse map.
$F$ is assumed to satisfy the canonical restrictions, so $F^{-1}$ is uniquely defined under the same
restrictions. Also recall that there is a uniquely defined Cremona inverse matrix $A_{F^{-1}}$ (Theorem~\ref{teorema_fundamental})
satisfying the canonical restrictions.
Moreover, both are automatically cohesive (see \cite[Lemma 4.1]{birational-linear}.

As explained in the previous section,  the graph corresponding to a monomial Cremona set in degree $2$
has either a unique circuit and this circuit has odd length $r\geq 3$ or else degenerates into a loop of length $r=1$.
For convenience we will call {\em root circuit} this basic unique subgraph in both cases.
Also, when talking of the log-matrix in this case we always mean the normal form $\mathcal{N}_F$ as established in
Lemma~\ref{format_logmatrix}.

\begin{Theorem}\label{degree_of_inverse}
Let $F=\{\xx^{v_1},\ldots, \xx^{v_n}\}\subset k[x_1,\ldots,x_n] (n\geq 2)$ denote a monomial Cremona set in degree $2$
satisfying the canonical restrictions and let $G_F$ stand for the corresponding graph.
Let $r$ denote the length of the root circuit of $G_F$
and let $s$ stand for the number of vertices
of $G_F$ off the root circuit whose vertex degree is $\geq 2$.
Then:
\begin{enumerate}
\item[{\rm (a)}] The degree of the Cremona inverse $F^{-1}$
is $(r+1)/2 + s\,${\rm ;}
\item[{\rm (b)}] The entries  of the Cremona inverse matrix $A_{F^{-1}}$ satisfy the following conditions:
    \begin{enumerate}
    \item[{\rm (i)}] The entries on the the main diagonal
are all nonzero
    \item[{\rm (ii)}] All entries belong to $\{0,1,2\}${\rm ;} moreover, $2$ is an entry on
    the $i$th row if and only if the corresponding vertex $x_i$ does not belong to the root
    circuit of $G_F$ and has degree $\geq 2${\rm ;} in particular, the inverse map is defined by
    squarefree monomials if and only if the root circuit of $G_F$ has empty second neighborhood.
    \end{enumerate}
\item[{\rm (c)}]  Write $\widetilde{G_{F}}$ for the subgraph obtained from $G_F$ by omitting the vertices
of degree $1$ and the corresponding incident edges and let $\{f_{i_1}, \ldots, f_{i_m}\}$ denote the monomials corresponding
to the edges of $\widetilde{G_{F}}$ not belonging to the root circuit $\{x_1,\ldots,x_r\}$ of $G_F$.
Then the inversion monomial factor is $x_1\cdots x_r f_{i_1} \cdots f_{i_m}$, where the order of the vertices of the circuit
and of the other edges is as in the normal form ${\cal N}_F$.
\end{enumerate}
\end{Theorem}
\demo
Let us first focus on the squarefree case (i.e., when $G_F$ is a simple graph). The degenerate case of a loop will need
just a few adjustments to be explained at the end.

Drawing on Theorem~\ref{cremona_grafo}, one considers the family of all connected graphs on $n$ vertices, having a unique circuit
of a fixed length $3\leq r\leq n$ and
proceed by induction on $n\geq r$.

For $n=r$, the graph reduces to a circuit of length $r$. The  Cremona inverse set satisfying the canonical restrictions is well-known:
it comes from the set of minimal coverings of the circuit by taking, for each such covering, the product of the
variables corresponding to the vertices
of the covering (see \cite[2.2, p. 347]{RuSi}).
Explicitly, one may take
\begin{equation}\label{coverings}
F^{-1}=\{y_1y_3\cdots y_n, \; y_2y_4\cdots y_{n-1}y_1,\;  \ldots, \; y_{n-1}y_ny_1\cdots y_{n-3},\; y_ny_0y_2\cdots y_{n-2}\}.
\end{equation}
Clearly, the degree is  $(n+1)/2$ and the entries on the main diagonal of $A_{F^{-1}}$ are all equal to $1$.
Further, the entries of the Cremona inverse matrix are $0,1$ and the monomial inversion factor
is $x_1\cdots x_n$ as is seen directly. This takes care of statements (a) through (c) in this
particular situation.

Now assume that $n\geq r+1$.
Then there exists  a vertex off the circuit whose degree in $G_F$ is $1$ - any vertex lying on the last
nonempty neighborhood will do.
Moreover, by suitably reordering the indices corresponding to the vertices on this neighborhood, one may assume that
the vertex corresponds to $x_n$ and the unique edge issuing from it corresponds to the monomial $\xx^{v_n}=x_jx_n$,
where $j\in\{1,\ldots, n-1\}$ is such that $x_j$ corresponds to the unique vertex adjacent to $x_n$ -- note that this
vertex belongs to the next to last neighborhood of the circuit.

Clearly, the set $F'=F\setminus \xx^{v_n}=\{\xx^{v_1},\ldots, \xx^{v_{n-1}}\}\subset k[\xx\setminus x_n]=k[x_1,\ldots,x_{n-1}]$
is still a Cremona set in degree $2$ satisfying the canonical restrictions.
By the inductive hypothesis, its Cremona inverse set $F'^{-1}$ and inversion vector $\gamma'$ satisfy statements (a) through (c).
In particular, the degree of $F'^{-1}$ is
$(r+1)/2+s'$, where $s'$ is the number of vertices
of $G_{F'}$ off the circuit whose vertex degree in $G_{F'}$ is $\geq 2$.

Since we have not messed with the normal arrangement of the earlier neighborhoods of the circuit,
 the normal form of the log-matrix of $F$ has the shape
$${\cal N}_F = \left(
\begin{array}{c@{\quad\vrule\quad}c} 
{\cal N}_{F'} & M \\
\multispan2\hrulefill\\
0 & 1 \\
\end{array}
\right),$$
where ${\cal N}_{F'}$ is a normal form of the log-matrix of $F'$ and $M$ denotes an $(n-1)\times 1$
matrix whose only nonzero entry is its $j$th entry.

We will now argue that the Cremona inverse matrix of ${\cal N}_F$ inherits a similar shape, with ${\cal N}_{F'}$ replaced by
the Cremona inverse matrix of ${\cal N}_{F'}$.
From this shape, we will read off the degree of $F^{-1}$ and the remaining assertions in the statements.

For this, one has to analyze the role of the special vertex corresponding to $x_j$, where $j\in\{1,\ldots,n-1\}$.
We accordingly divide this analysis into three distinct cases, to wit:

\medskip

\noindent{\bf Case 1.} The vertex corresponding to $x_j$ belongs to the common circuit of $G_F$ and $G_{F'}$.

\medskip

In this case, since $x_j$ belongs to the common circuit of $G_F$ and $G_{F'}$, one has $s'=s$ and,
clearly, $\widetilde{G_{F}}=\widetilde{G_{F'}}$ since $x_j$ has degree $\geq 2$.

Further, since $r\geq 3$, one can rewrite a normal form of $F'$ in which $j=2$ and $\xx^{v_1}=x_1x_2, \xx^{v_2}=x_2x_3$
are the two adjacent edges of the circuit.
The claim is that the following matrix has the properties in the statement of Theorem~\ref{teorema_fundamental} with
regard to ${\cal N}_F$, i.e., is the Cremona inverse matrix of ${\cal N}_F$.
$$B = \left(
\begin{array}{c@{\quad\vrule\quad}c} 
A_{F'^{-1}} & N \\
\multispan2\hrulefill\\
0 & 1 \\
\end{array}
\right)$$
where $A_{F'^{-1}}$ is the log-matrix of a uniquely defined inverse to $F'$ as in [loc.cit.], and $N$ is the $(n-1)\times 1$
matrix in which
$$N_l = \left\{
\begin{array}{ll} 
(A_{F'^{-1}})_{_{l1}} - 1, & \mbox{ if } l = 1\\
(A_{F'^{-1}})_{_{l1}}, & \mbox{ otherwise. }\\
\end{array}
\right.
$$
Note that subtracting $1$ in the first alternative above makes sense since, by the inductive hypothesis, no entry along the main
diagonal of $A_{F'^{-1}}$ is null.

Let ${\cal N}_{F'}\cdot A_{F'^{-1}}=\Gamma' +I_{n-1}$ be the fundamental matrix equation of inversion as in
Theorem~\ref{teorema_fundamental} relative to $F'$, with $\Gamma'=[\underbrace{\gamma'|\cdots |\gamma'}_{n-1}]$.
The usual block multiplication then yields
$${\cal N}_F\cdot B = \left(
\begin{array}{c@{\quad\vrule\quad}c} 
{\cal N}_{F'}\cdot A_{F'^{-1}} & {\cal N}_{F'}\cdot N+M \\
\multispan2\hrulefill\\
0 & 1 \\
\end{array}
\right) = \left(
\begin{array}{c@{\quad\vrule\quad}c@{\quad\vrule\quad}c@{\quad\vrule\quad}c} 
\gamma' & \ldots & \gamma' & {\cal N}_{F'}\cdot N+M\\
\multispan4\hrulefill\\
0 & \ldots & 0 & 0 \\
\end{array}
\right) + I_n,$$
where $M$ is $(n-1)\times 1$  and  $\gamma'=[(\gamma')_1,\cdots ,(\gamma')_{n-1}]^t$.
A straightforward calculation now yields $({\cal N}_{F'}\cdot N+M)_{_l}=(\gamma')_{_l}$, for $l=1,\ldots, n-1$. Indeed:

\begin{itemize}
\item If $l \neq 1,2$ one has $M_l = 0$ and  $({\cal N}_{F'})_{_{l1}}=0$, hence
$$({\cal N}_{F'}\cdot N+M)_l = \sum_{k=1}^{n-1} ({\cal N}_{F'})_{_{lk}}n_{_k} = \sum_{k=1}^{n-1}
({\cal N}_{F'})_{_{lk}}(A_{F'^{-1}})_{_{k1}} = (\gamma')_{_l}.$$

\item If $l=2$,  it is the case that $M_2=1$ and  $({\cal N}_{F'})_{_{21}} = 1$, therefore
$$({\cal N}_{F'}\cdot N+M)_{_2} = \sum_{k=1}^{n-1} ({\cal N}_{F'})_{_{2k}}n_{_k} + M_{_2} = \sum_{k=1}^{n-1} ({\cal N}_{F'})_{_{2k}}
(A_{F'^{-1}})_{_{k 1}} - ({\cal N}_{F'})_{_{21}} + M_{_2} = (\gamma')_{_2}.$$

\item If $l=1$ then $M_1=0$ and $({\cal N}_{F'})_{_{11}} = 1$ and hence
\begin{eqnarray*}({\cal N}_{F'}\cdot N+M)_{_{1}}& = &\sum_{k=1}^{n-1} ({\cal N}_{F'})_{_{1k}}n_{_k} =
\sum_{k=1}^{n-1} ({\cal N}_{F'})_{_{1k}}(A_{F'^{-1}})_{_{k1}} -
({\cal N}_{F'})_{_{11}} \\
&= &(\gamma')_{_{1}} + 1 - ({\cal N}_{F'})_{_{11}} = (\gamma')_{_1}.
\end{eqnarray*}
\end{itemize}

Hence this will yield an inversion equation ${\cal N}_F\cdot B=\Gamma +I_n$, with $\Gamma=[\gamma|\cdots |\gamma]$,
where $\gamma_l=(\gamma')_l$ for $l=1 ,\ldots,n-1$ and $\gamma_n=0$.
Moreover, since $F'^{-1}$ satisfies the canonical restrictions, so does $B$ by construction.
By Theorem~\ref{teorema_fundamental}, $B$ is the uniquely defined matrix giving the inverse to $F$, as stated.
Moreover, by the explicit format of $B$ with the nature of the entries of $N$, the inductive hypothesis
implies that the entries of $B$ are
$0,1,2$, while $2$ appears exactly on the $i$th row if and only if $x_i$ is a vertex off the root circuit
with degree $\geq 2$.
Since $|\gamma|=|\gamma'|$ and, besides, no entry on the main diagonal of $B$ is null, we are through in this case
for all statements of the theorem.

\medskip

\noindent{\bf Case 2.} The vertex corresponding to $x_j$ does not belong to the circuit and its degree in $G_{F'}$ is $\geq 2$

\medskip

In this case, since $G_{F'}$ only misses the vertex $x_n$ and its unique adjacent edge $x_jx_n$, and since $s$ only counts
vertices of degree $\geq 2$, it follows that $s'=s$ in this case too, and further,
$\widetilde{G_{F}}=\widetilde{G_{F'}}$ since $x_j$ has degree $\geq 2$.

By hypothesis, there are at least two indices $1\leq i' < j < i\leq n-1$ such that $x_{i'}x_j$ and $x_jx_i$ belong to the set $F$.
Note that the vertex $x_{i'}$ might belong to the circuit.
We claim that the following matrix has the properties in the statement of Theorem~\ref{teorema_fundamental} with
regard to ${\cal N}_F$, i.e., is its Cremona inverse matrix:
$$B = \left(
\begin{array}{c@{\quad\vrule\quad}c} 
A_{F'^{-1}} & N \\
\multispan2\hrulefill\\
0 & 1 \\
\end{array}
\right)$$
where $A_{F'^{-1}}$ is the Cremona inverse matrix of ${\cal N}_{F'}$ as in [loc.cit.], and $N$ is the $(n-1)\times 1$
matrix in which
$$N_l = \left\{
\begin{array}{ll} 
(A_{F'^{-1}})_{_{li}} - 1, & \mbox{ if } l = i\\
(A_{F'^{-1}})_{_{li}}, & \mbox{ otherwise. }\\
\end{array}
\right.
$$
where $i$ is the largest of the two indices taken above.
Let ${\cal N}_{F'}\cdot A_{F'^{-1}}=\Gamma' +I_{n-1}$ be the fundamental matrix equation of inversion as in
Theorem~\ref{teorema_fundamental} relative to $F'$, with $\Gamma'=[\underbrace{\gamma'|\cdots |\gamma'}_{n-1}]$.
Block multiplication then yields
$${\cal N}_F\cdot B = \left(
\begin{array}{c@{\quad\vrule\quad}c} 
{\cal N}_{F'}\cdot A_{F'^{-1}} & {\cal N}_{F'}\cdot N+M \\
\multispan2\hrulefill\\
0 & 1 \\
\end{array}
\right) = \left(
\begin{array}{c@{\quad\vrule\quad}c@{\quad\vrule\quad}c@{\quad\vrule\quad}c} 
\gamma' & \ldots & \gamma' & {\cal N}_{F'}\cdot N+M\\
\multispan4\hrulefill\\
0 & \ldots & 0 & 0 \\
\end{array}
\right) + I_n,$$
where $M$ is $(n-1)\times 1$  and $\gamma'=[(\gamma')_1,\cdots ,(\gamma')_{n-1}]^t$.
A straightforward calculation now yields $({\cal N}_{F'}\cdot N+M)_{_l}=(\gamma')_{_l}$, for $l=1,\ldots, n-1$:

\begin{itemize}
\item If $l \neq j,i$ one has $M_l = 0$ and  $({\cal N}_{F'})_{_{li}}=0$, hence
$$({\cal N}_{F'}\cdot N+M)_l = \sum_{k=1}^{n-1} ({\cal N}_{F'})_{_{lk}}n_{_k} = \sum_{k=1}^{n-1}
({\cal N}_{F'})_{_{lk}}(A_{F'^{-1}})_{_{ki}} = (\gamma')_{_l}.$$

\item If $l=j$,  it is the case that $M_j=1$ and  $({\cal N}_{F'})_{_{ji}} = 1$, therefore
$$({\cal N}_{F'}\cdot N+M)_{_j} = \sum_{k=1}^{n-1} ({\cal N}_{F'})_{_{jk}}n_{_k} + M_{_j} = \sum_{k=1}^{n-1} ({\cal N}_{F'})_{_{jk}}
(A_{F'^{-1}})_{_{k i}} - ({\cal N}_{F'})_{_{ji}} + M_{_j} = (\gamma')_{_j}.$$

\item If $l=i$ then $M_i=0$ and $({\cal N}_{F'})_{_{ii}} = 1$ and hence
\begin{eqnarray*}({\cal N}_{F'}\cdot N+M)_{_{i}}& = &\sum_{k=1}^{n-1} ({\cal N}_{F'})_{_{ik}}n_{_k} =
\sum_{k=1}^{n-1} ({\cal N}_{F'})_{_{ik}}(A_{F'^{-1}})_{_{ki}} -
({\cal N}_{F'})_{_{ii}} \\
&= &(\gamma')_{_{i}} + 1 - ({\cal N}_{F'})_{_{ii}} = (\gamma')_{_i}.
\end{eqnarray*}
\end{itemize}

Hence, this will yield an inversion equation ${\cal N}_F\cdot B=\Gamma +I_n$, with $\Gamma=[\gamma|\cdots |\gamma]$,
where $\gamma_l=(\gamma')_l$ for $l=1 ,\ldots,n-1$ and $\gamma_n=0$.
The conclusion is identical to the one in the previous case for all three statements
(a) through (c), so we are done in this case as well.

\medskip

\noindent{\bf Case 3.} The vertex corresponding to $x_j$ does not belong to the circuit and its degree in $G_{F'}$ is $1$

\medskip

In this case, one readily sees that $s=s'+1$ and, by a similar token, $\widetilde{G_{F}}=\widetilde{G_{F'}}\cup \{x_ix_j\}$,
where $x_i$ is the unique vertex
of $G_{F'}$ adjacent to $x_j$.
Note that, in this situation, we need to prove that the vector $\gamma$ that appears in the
inversion equation of $F$ has modulo $|\gamma'|+2$, where $\gamma'$ is the corresponding vector for the inversion equation of $F'$;
more precisely, we need $\xx^{\gamma} = \xx^{\gamma'}x_ix_j$.

Now, let again $i<j$ denote the unique index such that $x_ix_j$ belongs to $F'$ --
here $x_i$ may or may not belong to the circuit.
In this case we find it appropriate to express the Cremona inverse matrix of ${\cal N}_F$ in the form
$$B = \left(
\begin{array}{c@{\quad\vrule\quad}c} 
A_{F'^{-1}} & N \\
\multispan2\hrulefill\\
0 & 1 \\
\end{array}
\right) + E,$$
where $N$ is the $i$th column of $A_{F'^{-1}}$ and $E=(e_{kl})$ is the $n\times n$ matrix defined by
$$e_{kl} = \left\{
\begin{array}{ll} 
1, & \mbox{ if } k=j, l\neq n\\
0, & \mbox{ otherwise.}\\
\end{array}
\right.
$$
Multiplying we find
$${\cal N}_F\cdot B = \left(
\begin{array}{c@{\quad\vrule\quad}c} 
{\cal N}_{F'}\cdot A_{F'^{-1}} & {\cal N}_{F'}\cdot N+M \\
\multispan2\hrulefill\\
0 & 1 \\
\end{array}
\right) + \left(
\begin{array}{ccc@{\quad\vrule\quad}c} 
\alpha_j & \ldots & \alpha_j & \mathbf{0}\\
\multispan4\hrulefill\\
0 & \ldots & 0 & 0 \\
\end{array}
\right)= $$
$$= \left(
\begin{array}{c@{\quad\vrule\quad}c@{\quad\vrule\quad}c@{\quad\vrule\quad}c} 
\gamma' + \alpha_j & \ldots & \gamma' + \alpha_j & {\cal N}_{F'}\cdot N+M\\
\multispan4\hrulefill\\
0 & \ldots & 0 & 0 \\
\end{array}
\right) + I_n,$$
where $\alpha_j$ is the $j$th column of ${\cal N}_{F'}$.

We now assert that:
\begin{enumerate}
\item ${\cal N}_{F'}\cdot N+M = \gamma' + \alpha_j$
\item Every row of $B$ has a zero entry.
\end{enumerate}
 Once these are settled, the above gives an inversion equation
for $F$ with unique vector $\gamma:=\gamma' + \alpha_j$ and, since $\alpha_j$ is a column of a Cremona matrix in degree $2$, it
will follow that $\xx^{\gamma} = \xx^{\gamma' + \alpha_j}
= \xx^{\gamma'}x_ix_j$ -- in particular, $|\gamma|=|\gamma'|+2$ -- as was to be shown.
It is also clear from the form of $B$ that there is no null entry along its main diagonal.
Moreover, since $x_j$ has degree $1$ on $G_{F'}$ the entries on the $j$th row of $A_{F'^{-1}}$
are $0,1$ and not all are zero.
Once more, the explicit format of $B$ shows that its entries are $0,1,2$, while $2$ appears exactly on the $i$th row
if and only if $x_i$ is a vertex off the root circuit
with degree $\geq 2$.

To prove the first assertion, we proceed again along several cases:

\begin{itemize}
\item If $l \neq j,i$ one has $M_l = 0$ and  $({\cal N}_{F'})_{_{lj}}=0$, hence
$$({\cal N}_{F'}\cdot N+M)_l =  \sum_{k=1}^{n-1}
({\cal N}_{F'})_{_{lk}}(A_{F'^{-1}})_{_{ki}} = (\gamma')_{_l} +\alpha_{lj}.$$

\item If $l=j$,  it is the case that $M_j=1$ and  $({\cal N}_{F'})_{_{jj}} = 1$, therefore
$$({\cal N}_{F'}\cdot N+M)_{_j} = \sum_{k=1}^{n-1} ({\cal N}_{F'})_{_{jk}}
(A_{F'^{-1}})_{_{k i}} +1= (\gamma')_{_j} + 1= (\gamma')_{_j} + \alpha_{jj}.$$

\item If $l=i$ then $M_i=0$ and $({\cal N}_{F'})_{_{ij}} = 1$, hence
$$({\cal N}_{F'}\cdot N+M)_{_{i}} = \sum_{k=1}^{n-1} ({\cal N}_{F'})_{_{ik}}(A_{F'^{-1}})_{_{ki}}
= (\gamma')_{_{i}} + 1 = (\gamma')_{_i} + \alpha_{ij}.
$$
\end{itemize}

As for the second assertion, it is obvious for any row except possibly for the $j$th row.
For the latter, it suffices to show that $(A_{F^{-1}})_{jn}=(A_{F'^{-1}})_{ji}=0$.
First observe that $({\cal N}_{F'})_{jj}=1$ and $({\cal N}_{F'})_{jk}=0$ for every $k\in \{1,\ldots n-1\}\setminus \{j\}$
since we are in the normal form and $j$ is the index
of a vertex of degree $1$ off the circuit.
Then the inversion equation of $F'$ implies that $(A_{F'^{-1}})_{jk}= \gamma'_j + \delta_{jk}$
for $k\in \{1,\ldots n-1\}$.
Since $A_{F'^{-1}}$ satisfies the canonical restrictions, some entry along its $j$th row is null.
Therefore $(A_{F'^{-1}})_{jk}=0$ for all $k\neq j$ and $(A_{F'^{-1}})_{jj}=1$; in particular, $(A_{F'^{-1}})_{ji}=0$.

\bigskip

To conclude, we explain the adjustment in the case where the circuit degenerates into a loop.
This concerns only the initial step in the induction process.
Since $F$ satisfies the canonical restrictions, with  $n\geq 2$ by assumption,  its constituents have no proper common factor.
Therefore, $n\geq s_1+2$ where, we recall, $s_1\geq 1$ stands for the set of edges of $G_F$ in the
first neighborhood of the loop.
Thus, the initial step could be vacuous or, alternatively, would start from $n= s_1+2$, while in the inductive step one would
then assume that $n> s_1+2$.
Taking up the second alternative, the initial step has $F=\{x_1^2,x_1x_2,\ldots, x_1x_{n-1}, x_{n-1}x_n\}$, with $n\geq 3$.
Consider the following set in degree $2$:
\begin{equation}\label{inverse_loop}
\{x_1x_{n-1},x_2x_{n-1},\ldots, x_{n-2}x_{n-1}, x_{n-1}^2, x_1x_n\},
\end{equation}
 where the roles of
$x_1$ and $x_{n-1}$ have been interchanged and the loop has moved to another slot in the sequence.
A direct calculation show that this set is the Cremona inverse set of $F$ -- one can compose the two sets on the nose
or else pass to the respective log-matrices and multiply them out to get the
inversion equation, with inversion vector $\gamma=(2,0,\ldots,0, 1,0)^t$.

The rest of the argument stays unchanged, as far as the inductive step goes.
\qed

\subsection{Towards a classification of Cremona maps of degree $2$}

In this part we first briefly state the types of graphs corresponding to the classification
suggested in \cite[5.1.2]{birational-linear}.

\begin{Definition}\label{cremona_types}\rm
\begin{enumerate}
\item A set of squarefree monomials satisfying the canonical restrictions
is called {\it doubly-stochastic\/} if its log-matrix
is doubly-stochastic, i.e., the entries of each column sum up to an integer $d\geq 1$
(i.e., the monomials have fixed degree $d$) and so do the entries of each row
(i.e., no variable is privileged or, the ``incidence'' degree of any variable is
also $d$).
\item A Cremona set satisfying the canonical restrictions is called a {\it $p$-involution\/} if it coincides with its
inverse set up to permutation on the source and the target.
\item A Cremona set satisfying the canonical restrictions is called {\it
apocryphal\/} if its inverse set has at least one non-squarefree monomial.
\end{enumerate}
\end{Definition}
The notion of a $p$-involution  has been introduced
in parallel to the classical situation of an involuting  Cremona map (up to a projective change of coordinates).
Here $p$ stands as short reminder for ``permutation''.
Note it makes perfect sense even if the set contains monomials which are not squarefree.

The next result tell us about the nature of a $p$-involution in degree $2$.

\begin{Proposition}\label{involution_deg_2}
Let $F$ denote a Cremona set of degree $2$ and let $G_F$ for the corresponding
graph.
Let $r$ denote the size of the unique circuit  {\rm (}possibly a loop{\rm )} in $G_F$
 and let $s$ stand for the number of vertices
of $G_F$ off the circuit whose vertex degree is $\geq 2$.
The following conditions are equivalent:
\begin{enumerate}
\item[{\rm (a)}] $F$ is a $p$-involution
\item[{\rm (b)}] The inverse $F^{-1}$ has degree $2$
\item[{\rm (c)}] Either $r=3$ {\rm (}triangle{\rm )} and $s_2=0$,
or else $r=1$  {\rm (}loop{\rm )} and $s=1$.
\end{enumerate}
\end{Proposition}
\demo The implication (a) $\Rightarrow$ (b) is trivial.

The implication (b) $\Rightarrow$ (c) is obtained as follows.
By Theorem~\ref{degree_of_inverse}, $(r+1)/2+s=2$.
Therefore, either $r=3$ and $s=0$ -- this corresponds to a circuit of length $3$ and possibly
additional edges all adjacent to the circuit, -- or else $r=1$ and $s=1$. The latter case means that the
corresponding graph consists of a loop and at least one edge adjacent
 to the loop and, moreover, only one of these edges has adjacent edges in the second neighborhood of the loop.

To see that (c) implies (a) we separate the two cases.

First take the case where $G_F$ have a circuit $C_3$ of length $3$.
By assumption, the neighborhood of order $2$ of $C_3$ is empty.
Therefore, the normal form of the log-matrix of $F$ has at most $2$ blocks:
$${\cal N}_F = \left(
\begin{array}{c@{\quad\vrule\quad}c}
{\cal N}_{C_3}\!\!\!\! & M \\
\multispan2\hrulefill\\
0  & I_t \\
\end{array}
\right),$$
where $t$ is the cardinality of the first neighboorhood of $C_3$ and $M$ is a $3\times t$ matrix having
exactly one nonzero entry on every column, this entry being $1$.

Moreover, the Cremona defined by $C_3$ is a $p$-involution by direct inspection (or as a trivial case
of (\ref{coverings})) and the inversion equation is
$${\cal N}_{C_3}\cdot {\cal N}_{C_3}\,'=[\gamma_{_{C_3}}|\gamma_{_{C_3}} |\gamma_{_{C_3}}]+I_3,$$
where $\gamma_{_{C_3}}=(1,1,1)^t$ and ${\cal N}_{C_3}\,'$ is obtained from ${\cal N}_{C_3}$
by applying the permutation $1\mapsto 2\mapsto 3\mapsto 1$ to its columns.

Consider the matrix
$$B= \left(
\begin{array}{c@{\quad\vrule\quad}c}
{\cal N}_{C_3}\,'\!\!\!\! & N \\
\multispan2\hrulefill\\
0  & I_t \\
\end{array}
\right),$$
where $N$ is to be determined so that there is an equality
${\cal N}_F\cdot B= [\gamma|\ldots |\gamma]+I_t$, with $\gamma=(1,1,1,0,\ldots,0)^t$.
Note that
$${\cal N}_F\cdot B=
\left(
\begin{array}{c@{\quad\vrule\quad}c} 
[\gamma_{_{C_3}}|\gamma_{_{C_3}}|\gamma_{_{C_3}}]+I_3 & {\cal N}_{C_3}\cdot N+M \\
\multispan2\hrulefill\\
0 & I_t \\
\end{array}
\right).
$$
That is, we are  to solve the equation
${\cal N}_{C_3}\cdot N+M=[\,\mathbb{I}|\ldots|\mathbb{I}\,]$ for $N$,
where $\mathbb{I}=(1,1,1)^t$.

Now, let $(\mathfrak{n}_{1j},\mathfrak{n}_{2j},\mathfrak{n}_{3j})^t$ denote the $j$th column of $N$.
Then the $j$th column of ${\cal N}_{C_3}\cdot N$ is $(\mathfrak{n}_{1j}+\mathfrak{n}_{3j},\,\mathfrak{n}_{1j}+\mathfrak{n}_{2j},\,
\mathfrak{n}_{2j}+\mathfrak{n}_{3j})^t$.
Since every column of $M$ has exactly one nonzero entry, and this entry is $1$, we are typically led to solve the system of equations
$$\left\{
\begin{array}{c}
\mathfrak{n}_{1j}+\mathfrak{n}_{3j}+1=1\\
\mathfrak{n}_{1j}+\mathfrak{n}_{2j}=1\\
\mathfrak{n}_{2j}+\mathfrak{n}_{3j}=1
\end{array}
\right.
$$
The solution is immediately seen to be $\mathfrak{n}_{1j}=\mathfrak{n}_{3j}=0,\,\mathfrak{n}_{2j}=1$.
Thus, $N$ is uniquely obtained and, like $M$, it has exactly one nonzero entry and this entry is $1$.
Moreover, by an obvious symmetry of the solution, applying the permutation $1\mapsto 2\mapsto 3\mapsto 1$
this time around to the rows of $N$, we see that the columns of $M$ and $N$ are the same.
This proves that $F$ is a $p$-involution and $B$ is the corresponding Cremona inverse matrix.

The loop case has already been described in the proof of  Theorem~\ref{degree_of_inverse} (see (\ref{inverse_loop})).
\qed

\bigskip

Now specialize to {\em squarefree} Cremona sets of degree $2$.
We assume throughout that the Cremona set in degree $2$ satisfies the canonical restrictions.
\begin{Definition}\label{short_and_general}\rm
A squarefree Cremona set in degree $2$ is of {\em short type} if
the (odd) circuit of the corresponding graph has empty second neighborhood; otherwise we say that the set
is of {\em long type} or of {\em general type}.
\end{Definition}

By Theorem~\ref{degree_of_inverse}, the Cremona set is of short type
if and only if the degree of its inverse is $(r+1)/2$, where $r$ is the length of the unique (odd) circuit
in the corresponding simple graph.

We now file a couple of consequences.

\begin{Corollary}\label{graph_types}
Let $F$ be a squarefree Cremona set in degree $2$ of short type
which is not doubly-stochastic.
Then the corresponding graph $G_F$ is a circuit of length $\geq 3$ with nonempty first
neighborhood. Moreover, $F$ is a $p$-involution if and only if the circuit has length $3$.
\end{Corollary}
\demo
It is evident that $G_F$ has the stated form.
The characterization of a $p$-involution is the content of the equivalence (a) $\Leftrightarrow$ (c) in Proposition~\ref{involution_deg_2},
\qed

\medskip

The nex result follows immediately from Theorem~\ref{degree_of_inverse} (b)(ii), but we wish to isolate it as
natural complement to the previous proposition.

\begin{Corollary}\label{graph_of_apocryphal}
Let $F$ be a squarefree Cremona set in degree $2$.
Then $F$ is apocryphal if and only if it is of general type.
\end{Corollary}

\begin{Remark}\rm It would be interesting to know if there is a ``fractalization'' of the different subtypes of a degree $2$
Cremona set of general type. The results so far seem to point in the direction that all ``look alike''.
\end{Remark}

\section{The role of Hilbert bases}

In this section we answer a question posed by R. Villarreal (oral communication) about the connection
between monomial Cremona transformations and Hilbert bases.

For the reader's convenience we review the needed background on this combinatorial topic, our main references being
\cite{BrGu}, \cite{Schr}, \cite{VillaBook} and \cite{comica}.

\subsection{Review of main facts}

We assume the elementary notions of polyhedral combinatorics.

Recall the partial order on $\mathbb{R}^n$ defined by $a = (a_1, \ldots, a_n)\leq c = (c_1, \ldots, c_n)$ if
$a_i \leq c_i$ for every $i$.
Given vectors $a,b\in \mathbb{R}^n$, their inner product will be denotes by $<a,b>$.
A hyperplane $H = H(a,c) = \{x \in \mathbb{R}^n| \textless x,a \textgreater = c\}\subset \mathbb{R}^n$ determines
two closed half-spaces
$$H^{+}(a,c) = \{x \in \mathbb{R}^n| \textless x,a \textgreater \geq c\} \mbox{ and } H^{-}(a,c) = \{x \in \mathbb{R}^n|
\textless x,a \textgreater \leq c\}.$$
We follow common usage of writing $H^+_a=H^+(a,0)$ and $H^-_a=H^-(a,0)$ when the hyperplane goes through the origin.

A ({\em convex}) {\em cone} in $\mathbb{R}^n$ is a nonempty set $C\subset \mathbb{R}^n$ such that, for all $x,y \in C$
and all real $\lambda, \mu \geq 0$, one has $\lambda x + \mu y \in C$.

A cone $C$ is {\em polyhedral} if it there exists a finite subset $S=\{s_1,\ldots,s_q\}\subset C$ such that
$C=\R_+\,S:=\{\lambda_1s_1+\cdots +\lambda_qs_q,\, \lambda_j\in \R_+\}$.
We then refer to $S$ as a generating set of $C$.

The {\em dual} to a polyhedral cone $C\subset \R^n$ is defined as
$$C^*:=\{u\in {\R^n}^*=\mbox{\rm Hom}_{\R}(\R^n,\R)\,|\, <v,u>\leq 0, \, \forall v\in C\},$$
where $<\,,\,>$ denotes the ordinary pairing on $\R^n\times {\R^n}^*$.
It is immediate to see that it suffices to take the pairing over a finite set of cone generators of $C$.
Common practice  identifies $\R^n$ with its dual space by identifying a vector basis with its dual basis.
In this way, the pairing can be seen as the usual inner product and the dual to a cone can be considered in the same space.
As such, one has a representation $C^*=\cap_{s\in S}H_s^-$ as intersection of half-spaces through the origin.
Moreover, a fundamental result  going back to Farkas, Minkowski and Weyl (see \cite[Theorem 1.1.31]{comica}) says that a cone is polyhedral
if and only if it is the intersection of a finite set of closed half-spaces through the origin.
From this follows that the dual to a polyhedral cone is also a polyhedral cone and $(C^*)^*=C$ for every polyhedral cone.

This well-known dichotomy of representing a polyhedral cone, both as the intersection of closed half-spaces through the origin and
as the set of nonnegative linear combinations of a finite set of vectors is very useful.
The first of these representation allows to write a polyhedral cone in the form $\{v\in \R^n| Av\leq 0\}$, for some real matrix $A$.

In this vein, we say that $C$ is {\em pointed} if the the linear system $Ax=0$ has only the trivial solution $x=0$, i.e., $A$ has maximal rank.
In other words, a pointed polyhedral cone contains no straight lines.
Since this condition also means that $C\cap (-C)=\{0\}$, where $C$ has the obvious meaning, a pointed cone is also called {\em
strongly convex}.
To free ourselves from the matrix representation in the notion, we can use the characterization in \cite[Proposition 1.1.56]{comica}
to the effect that a polyhedral cone is pointed if and only if its dual has maximal dimension (i.e., $n$).

\subsection{Cremona maps out of Hilbert bases}

We now come to the main concept of this part.
As a matter of further notation, given a subset $\mathcal{A}\subset \mathbb{R}^n$, denote by $\Z\mathcal{A}$ (respectively,
$\NN\mathcal{A}$) the integer lattice generated by $\mathcal{A}$ (respectively, the set of  lattice elements with nonegative coefficients).

\begin{Definition}\label{hilbert_base}\rm
A finite subset $H\subset \mathbb{R}^n$ is a {\em Hilbert base} if $\,\mathbb{Z}^n \cap \mathbb{R}_+H = \mathbb{N}H.$
A polyhedral cone is said to admit a Hilbert base if it contains a Hilbert base and is generated by it.
\end{Definition}

Note that the definition implies that a Hilbert base is contained in $\Z^n$.

The fundamental results regarding Hilbert bases are as follows.
The first tells us that Hilbert bases are pretty ubiquitous and often uniquely defined.

\begin{Theorem}{\rm (\cite[Theorem 16.4]{Schr})}
A rational polyhedral cone $C$ admits a  Hilbert base.
If, moreover, $C$ is pointed then it contains a unique minimal such base
in the sense that no proper subset is a generating Hilbert base.
\end{Theorem}

The second result even gives a hint as to the nature of a minimal such base.

\begin{Theorem}{\rm (\cite[proof of Theorem $1$]{GerSeb})}\label{bh}
Let $H$ denote a Hilbert base of a pointed polyhedral cone and let $r$ stand for the rank
of the lattice generated by $H$.
Then $H$ has a subset of $r$ linearly independent vectors forming a Hilbert base.
\end{Theorem}

Recall that an ideal $I$ of a ring $R$ is {\em normal} if all its powers are integrally closed in $R$.
Let $k$ be an arbitrary field.
If $v=(a_1,\ldots a_n)\in \NN^n$ and $\XX=\{X_1,\ldots,X_n\}$ are indeterminates over $k$ then
we set $\XX^v:=X_1^{a_1}\cdots X_n^{a_n}\in k[\XX]$ for the associated monomial as introduced in the first section.

Our first main result of this section is the following.

\begin{Theorem}\label{main_normal}
Let $v_1, \ldots, v_q \in \mathbb{N}^n$ $(q\geq n)$ be given such that the associated monomials
$\XX^{v_1},\ldots, \XX^{v_q}$ have the same degree $d\geq 1$.
If the ideal $(\XX^{v_1},\ldots, \XX^{v_q})\subset k[\XX]$ is normal then
$\{(v_1,1), \ldots, (v_q,1)\}$ is a Hilbert base.
\end{Theorem}
\demo
We first note that, since the ring $k[\XX]$ is normal, the ideal $(\XX^{v_1},\ldots, \XX^{v_q})\subset k[\XX]$
is normal if and only if the Rees algebra of this ideal is
normal.
But the Rees algebra is isomorphic to the semigroup ring  $k[\XX, \,\XX^{v_1}T,\ldots, \XX^{v_q}T]\subset k[\XX,T]$.
Applying the well-known combinatorial criterion (see, e.g., \cite[Corollary 7.2.29]{{VillaBook}}), we find that this algebra is normal
if and only if
$$\mathbb{Z}H' \cap \mathbb{R}_+H' = \mathbb{N}H',$$
where $H'=\{e_1, \ldots, e_n\}\cup H$, $H=\{ (v_1,1), \ldots, (v_q,1)\} \subset \mathbb{N}^{n+1}$, and $\{e_1, \ldots, e_n,e_{n+1}\}$
stands for the canonical basis of $\Z^{n+1}$.

On the other hand, it is well-known that, as a consequence of the so-called Farkas Lemma (\cite[Corollary 1.1.29]{comica}, also
\cite[Corollary 7.1d]{Schr}), for any subset $\mathcal{A}\subset \Z^m$ the natural inclusions
$\Z\mathcal{A} \cap \Q_+\mathcal{A} \subset \Z\mathcal{A} \cap \R_+\mathcal{A}$ and
$\Z^m \cap \Q_+\mathcal{A} \subset \Z^m \cap \R_+\mathcal{A}$ are equalities.

We are therefore to prove that
$$\mathbb{Z}H' \cap \mathbb{Q}_+H' \subset \mathbb{N}H' \,\Rightarrow\, \mathbb{Z}^{n+1} \cap \mathbb{Q}_+H \subset \mathbb{N}H.$$

Now, $\mathbb{Z}H' = \mathbb{Z}^{n+1}$ as, e.g., $e_{n+1} = (v_1,1) - v_{1,1}e_1 + \cdots + v_{n,1}e_n$,
where $v_1=(v_{1,1},...,v_{n,1})$.

Thus, let $z \in \mathbb{Z}^{n+1} \cap \mathbb{Q}_+H$, say, $z = (z_1, \ldots,z_{n+1}) =
\lambda_1(v_1,1) + \ldots + \lambda_q(v_q,1)$.
Then
$$z_{n+1} = \sum_{i=1}^{q} \lambda_i \;\mbox{ and } |z| = (\sum_{i=1}^{q} \lambda_i)(d+1) = z_{n+1} (d+1).$$
Since $H\subset H'$, also $z \in \mathbb{Z}^{n+1} \cap \mathbb{Q}_+H'$, hence by the assumption
and the above remark, one can write
$$z= \alpha_1(v_1,1)+ \ldots + \alpha_q(v_q,1) + \beta_1 (e_1,0) + \ldots + \beta_n (e_n,0),$$
for suitable $\alpha_1, \ldots, \alpha_q, \beta_1, \ldots,\beta_n \in \mathbb{N}$.
But then $z_{n+1} = \sum_{i=1}^{q}\alpha_i$ and
$$|z|= (\sum_{i=1}^{q}\alpha_i)(d+1) + \sum_{j=1}^{n}\beta_j = z_{n+1}(d+1) + \sum_{j=1}^{n}\beta_j,$$
hence $\sum_{j=1}^{n}\beta_j =0$, that is, $\beta_j = 0$ for every $j$. Consequently, $z \in  \mathbb{N}H$
as was required to show.
\qed

\medskip

The following result is pretty elementary, but we give a proof for the sake of completeness.

\begin{Lemma}\label{ones_is_pointed}
Let $H=\{v_1,\ldots,v_q\}\subset \Z^n$ be an arbitrary subset.
Then $\{(v_1,1),\ldots, (v_q,1)\}\subset \Z^{n+1}$  generates a pointed cone.
If, moreover,  the associated monomials of the vectors in $H$ have the same degree then this cone has the same dimension as the
cone generated by $H$.
\end{Lemma}
\demo
To prove that the polyhedral cone generated by $(H,1):=\{(v_1,1),\ldots, (v_q,1)\}$ is pointed
one argues that it is strongly convex, namely, let
$$\sum_j a_j (v_j,1)=\sum_j b_j (v_j,1),$$
with $a_j\geq 0, \,\forall j$ and $b_j\leq 0,\, \forall j$.
Looking at the last coordinate, we get $\sum_j a_j=\sum_j b_j$.
Forcefully, $a_j= 0, \,\forall j$.

Assuming now that the matrix whose columns are the vectors in $H$ is $d$-stochastic for some $d\geq 1$,
the matrix
$$\begin{pmatrix}
v_1 & \cdots & v_n\\
d & \cdots & d
\end{pmatrix}
$$
the last row is the sum of the rows of $[v_1|\ldots |v_q]$, hence both have the same rank.
But this is also the rank over $\Q$ (hence, over $\Z$) of the matrix whose columns are $(v_1,1),\ldots, (v_q,1)$.

Thus, we are through.
\qed

\begin{Remark}\rm The second statement of the above lemma is actually \cite[Exercise 6.2.23]{comica}.
As to pointedness, one notes that, more generally, a finite set of vectors in $\Z^n$ generate
a pointed polyhedral cone if, for some $i\in\{1,\ldots,n\}$, their $i$th coordinates are positive natural numbers.
\end{Remark}

The main combinatorial result of this part now follows.

\begin{Theorem}\label{short_Hilbert_base}
Let $H=\{v_1,\ldots,v_q\}\subset \NN^n$ be such that the associated monomials have the same degree $d\geq 1$
and $\Z H$ has rank $n$.
If $(H,1)=\{(v_1,1),\ldots, (v_q,1)\}\subset \Z^{n+1}$ is a Hilbert base then there exists an $n\times n$ submatrix of
$[v_1|\ldots |v_q]$ whose determinant is $d$.
\end{Theorem}
\demo By Lemma~\ref{ones_is_pointed}, the cone generated by $(H,1)$ is pointed.
Therefore, by Theorem~\ref{bh}, $(H,1)$ admits a subset $(H',1)$ of $n$ linearly independent vectors which is
a Hilbert base.

We claim that the matrix whose columns are the vectors in $H'$ gives the required result.
For this it suffices to show an isomorphism of $\Z$-modules $\Z^n/\Z H'\simeq \Z/d\Z$.
First note that $\Z^n/\Z H'$ is a torsion $\Z$-module since $\Z H'$ has rank $n$.
On the other hand, one has an exact sequence of $\Z$-modules
$$0 \rightarrow T_{\Z}(\mathbb{Z}^{n+1} / \mathbb{Z} H) \rightarrow T(\mathbb{Z}^{n} / \mathbb{Z} H') \rightarrow
\mathbb{Z}/d\Z \rightarrow 0,$$
where $T_{\Z}$ denotes $\Z$-torsion (see, e.g., \cite[the proof of Theorem 1.1]{SiVi}).

Now, $T_{\Z}(\Z^{n+1}/\Z H)= \R H \cap \Z^{n+1}/\Z H$ by \cite[Lemma 1.2.11]{comica}.
Since $H$ is a Hilbert base, $\R_+ H \cap \Z^{n+1}=\NN H$.
But the latter equality implies the containment
$\R H \cap \Z^{n+1}\subset\Z H$, a fact that is readily checked by writing every
coefficient $a\in \R$ in the form $\lfloor a\rfloor +b$,
with $b\in \R_+$.
\qed

\begin{Corollary}\label{Hilbert2Cremona}
Let $H:=\{v_1,\ldots,v_q\}\subset \NN^n$ be such that the associated monomials have the same degree $d\geq 1$.
If $\Z H$ has rank $n$ and
$\{(v_1,1),\ldots, (v_q,1)\}\subset \Z^{n+1}$ is a Hilbert base then there exist $n$
vectors in $H$ such that the associated monomials
define a Cremona transformation of $\pp^{n-1}$.
\end{Corollary}
\demo By Theorem~\ref{short_Hilbert_base}, there exists an $n\times n$ submatrix of
$[v_1|\ldots |v_q]$ whose determinant is $d$.
Since $d\neq 0$, the respective associated monomials cannot have a proper common factor.
Therefore, the criterion of \cite[Lemma 2.2]{birational-linear} applies.
\qed

\medskip

The following result bundles up the previous results, its contents bridging between
combinatorics and birational geometry.

\begin{Theorem}\label{normal2Hilbert2Cremona}
Let $\XX^{v_1},\ldots, \XX^{v_q}\subset k[\XX]=k[X_1,\ldots,X_n]$ $(q\geq n)$ be monomials
of the same degree generating a normal ideal. Then there exist $n$
among these monomials defining a Cremona transformation of $\pp^{n-1}$.
\end{Theorem}
\demo
It follows immediately from the previous corollary and Theorem~\ref{main_normal}
\qed

\begin{Example}\rm
($n=3$) Consider the set of all monomials of degree $2$.
Clearly, these generate a normal ideal of $k[x,y,z]$.
Up to a permutation of variables and generators, there are
exactly two subsets of $3$ monomials defining a plane Cremona map
each. This is because these must coincide with the two known
quadratic Cremona maps with $3$ distinct base points and with $2$
distinct base points plus an infinitely near one, respectively:
$xy,xz,yz$ and $x^2,xy,yz$.
Both generate normal ideals, as is well-known or easy to check.
Thus, by Theorem~\ref{main_normal}, the corresponding sets of $\NN^4$
obtained by adding $1$ as the $4$th coordinate are Hilbert bases.
\end{Example}

\bibliographystyle{plain}

\bigskip

\noindent{\bf Addresses:}

\medskip

{\sc B. Costa}, Departamento de Matem\'atica, CCEN, Universidade
Federal de Pernambuco,  Pernambuco,
Brazil

{\it Email}: {\sc bcs@dmat.ufpe.br}

\bigskip

{\sc A. Simis}, Departamento de Matem\'atica, CCEN, Universidade
Federal de Pernambuco,  Pernambuco,
Brazil

{\it Email}: {\sc aron@dmat.ufpe.br}

\end{document}